\documentstyle[amscd,amssymb,amsthm,verbatim,amsopn,amsopn,12pt]{amsart}

\theoremstyle{plain}
\newtheorem{Thm}{Theorem}[section]

\newtheorem{Prop}[Thm]{Proposition}
\newtheorem{Cor}[Thm]{Corollary}
\theoremstyle{definition}

\newtheorem{Def}[Thm]{Definition}
\newcommand{\bnum}{\begin{enumerate}}
\newcommand{\enum}{\end{enumerate}}

\newcommand{\Z}{\mathbb{Z}} 
\newcommand{\N}{\mathbb{N}}

\begin{document}
\title[Generalized non-commuting graph of a finite ring]%
	{On generalized non-commuting graph of a finite ring}
\author[J. Dutta, D. K. Basnet and R. K. Nath]%
	{Jutirekha Dutta, Dhiren K. Basnet \& Rajat K. Nath*}
\thanks{Corresponding author}
\date{}
\maketitle
\begin{center}\small{\it Department of Mathematical Sciences, Tezpur University,}
\end{center}
\begin{center}\small{\it Napaam-784028, Sonitpur, Assam India.}
\end{center}

\begin{center}\small{\it Email: jutirekhadutta@@yahoo.com,  dbasnet@@tezu.ernet.in and rajatkantinath@@yahoo.com*}
\end{center}

\medskip

\noindent \textit{\small{\textbf{Abstract:}
 Let $S, K$ be two subrings of a finite ring $R$. Then the generalized non-commuting graph of subrings $S, K$ of $R$, denoted by $\Gamma_{S, K}$, is a simple graph whose vertex set is $(S \cup K) \setminus (C_K(S) \cup C_S(K))$ and two distinct vertices $a, b$ are adjacent if and only if $a \in S$ or $b \in S$ and $ab \neq ba$.   We  determine the diameter, girth and some dominating  sets for $\Gamma_{S, K}$. Some   connections between the $\Gamma_{S, K}$ and $\Pr(S, K)$ are also obtained. Further, $\Z$-isoclinism between two pairs of finite rings is defined and showed that  the generalized non-commuting graphs of two $\Z$-isoclinic pairs are isomorphic under some condition.  }}

\bigskip

 

\medskip

\section{Introduction}
Throughout this paper  $R$ denotes a finite ring, $S$ and  $K$ denote two subrings of $R$. Let $C_K(S) = \{k \in K : ks = sk \;\forall\; s \in S\}$ , $C_S(K) = \{s \in S : ks = sk \;\forall\; k \in K\}$. Note that $C_K(S)$ and $C_S(K)$ are subrings of $K$ and $S$ respectively. In this paper, we consider the graph $\Gamma_{S, K}$ associated to the subrings $S$ and $K$ of $R$ as follows: We take $(S \cup K) \setminus (C_K(S) \cup C_S(K))$ as the vertex set of $\Gamma_{S, K}$ and two distinct vertices $a$ and $b$ are adjacent if and only if $a \in S$ or $b \in S$ and $ab \neq ba$.

  It is clear that for any subrings $S, K$ of $R$ such that $S \subseteq K$, the vertex set of $\Gamma_{S, K}$ is $K \setminus C_K(S)$. Further, if $S = R$ then $\Gamma_{S, K}$ becomes $\Gamma_{R, R} = \Gamma_R$, the non-commuting graph of $R$.  The notion of non-commuting graph of a finite ring  was introduced by Erfanian, Khashyarmanesh and Nafar \cite{ekn15}. Many mathematicians have studied algebraic structures by means of graph theoretical properties in the last decades (see \cite{abdollahi07, abdollahi08, abdollahiakbarimaimani06, beck88, dn16, ov11} etc.). The motivation of this paper lies in the works of Erfanian et al. \cite{GET14,TE13}. 
 Also the techniques adopted to prove various results in this paper are more or less similar in nature  to that in \cite{TE13,GET14}.

 We recall that the commuting probability of a finite ring $R$ is the probability that a randomly chosen pair of elements of $R$ commute. That is,
\[
\Pr(R) = \frac{|\{(r, t) \in R \times R : rt = tr\}|}{|R||R|}. 
\]
This ratio was introduced by MacHale \cite{machale} and studied by MacHale et al. in \cite{BM, BMS}.  We generalize $\Pr(R)$    by the following ratio 
\[
\Pr(S, K) = \frac{|\{(s, k) \in S \times K : sk = ks\}|}{|S||K|}
\]
where $S$ and $K$ are two subrings of $R$. Various properties of $\Pr(S, K)$ are studied in \cite{ddn}. Clearly, $\Pr(R, R) = \Pr(R)$.  It may be mentioned here that $\Pr(S, K)$ when $K = R$ is studied in  \cite{duttabasnetnath}.

In Section $2$, we give some preliminary results regarding $\Gamma_{S, K}$. In Section $3$, we determine diameter, girth and some dominating sets for $\Gamma_{S, K}$.
In Section $4$, we derive some connections between $\Gamma_{S, K}$ and $\Pr(S, K)$. In the last section, we define $\Z$-isoclinism between two pairs of rings and find some connections between two isoclinic pairs of rings and their generalized non-commuting graphs. 

\section{Preliminary results}
In this section, we derive some preliminary results some of which are used in the forthcoming sections. For a graph $G$, we write $V(G)$ and $E(G)$ to denote the set of vertices    and set of edges   of $G$ respectively. We write $deg(v)$ to denote  the degree of the vertex $v$, which is the number of edges incident on $v$.

\begin{Prop}
Let $R$ be a non-commutative ring and $S, K$ two subrings of $R$. Let $r \in V(\Gamma_{S, K})$ then
\begin{enumerate}
\item $deg(r) = |S \cup K| - |C_S(r) \cup C_K(r) \cup C_S(K)|$ if $r \in S \setminus (K \cup C_K(S) \cup C_S(K)) = S \setminus (K \cup C_S(K))$.
\item $deg(r) = |S \cup K| - |C_S(r) \cup C_K(r)|$ if $r \in S \cap K$.
\item $deg(r) = |S| - |C_S(r) \cup C_{S \cap K}(S)|$ if $r \in K \setminus (S \cup C_K(S) \cup C_S(K)) = K \setminus (S \cup C_K(S))$.
\end{enumerate} 
\end{Prop}
\begin{pf}
The proof follows from the definition of $\Gamma_{S, K}$.
\end{pf}

As a consequence of the above proposition, we have the following corollary.
\begin{Cor}\label{deg_cor}
Let $R$ be a non-commutative ring with subrings $S$ and $K$ such that $S \subseteq K$. Then
\begin{enumerate}
\item $deg(r) = |K| - |C_K(r)|$ if $r \in V(\Gamma_{S, K}) \cap S$.
\item $deg(r) = |S| - |C_S(r)|$ if $r \in V(\Gamma_{S, K}) \cap (K \setminus S)$.
\item There is no isolated vertex in $V(\Gamma_{S, K})$.
\item $\Gamma_{S, K}$ is empty graph if and only if $S$ is commutative. 
\end{enumerate} 
\end{Cor}

Recall that a star graph is a tree on $n$ vertices in which one vertex has degree $n - 1$ and the others have degree $1$. A bipartite graph is a graph whose vertex set can be partitioned into two disjoint parts in such a way that the two end vertices of every edge lie in different parts. A  complete bipartite graph is a bipartite graph such that two vertices are adjacent if and only if they lie in different parts. In the following theorems we shall show that if $G$ is a star graph or complete bipartite graph or an $n$-regular graph, where $n$ is a square free odd positive integer, then $G$ can not be realized by $\Gamma_{S, K}$ for any two subrings $S, K$ of a ring $R$ such that $S\subseteq K$.

\begin{Thm}\label{star}
There is no non-commutative  ring $R$ with subrings $S, K$ and $S \subseteq K$ such that  $\Gamma_{S, K}$ is a star graph.
\end{Thm}
\begin{pf}
Suppose there exists a ring $R$ with non-commutative subrings $S, K$  such that $S \subseteq K$ and $\Gamma_{S, K}$ is a star graph. Then there exists a unique vertex of degree $|V(\Gamma_{S, K})| - 1$. Suppose first that $s \in S$ is that vertex. So $deg(s) = |V(\Gamma_{S, K})| - 1$ which gives $|C_K(S)|(|C_K(s)|/|C_K(S)| - 1) = 1$. So, $|C_K(S)| = 1$ and $(|C_K(s)|/|C_K(S)| - 1) = 1$. This gives $|C_K(s)| = 2$ and  $deg(r) = 1$ for any $r \in K \setminus S$. Therefore, $[S : C_S(r)] = |S|/(|S| - 1)$, which is not possible.

Next we suppose that $r \in K \setminus S$ is the unique vertex having degree $|V(\Gamma_{S, K})| - 1$. Then for any $s \in S \cap V(\Gamma_{S, K})$ we have $deg(s) = 1$ which again gives $[K : C_K(s)] = |K|/(|K| - 1)$, a contradiction.
This completes the proof.  
\end{pf}

\begin{Thm}\label{complete relative non-commuting}
There is no non-commutative ring $R$ with subrings $S, K$ and $S \subseteq K$ such that $\Gamma_{S, K}$ is complete bipartite.
\end{Thm}

\begin{pf}
Let $R$ be a finite non-commutative ring and $S, K$ two subrings of $R$, where  $S \subseteq K$, such that    $\Gamma_{S, K}$ is complete bipartite.
Then we have two disjoint subsets   $V_1$ and $V_2$ of $V(\Gamma_{S, K})$ such that $|V_1| + |V_2| = |K| - |C_K(S)|$. Suppose $S \cap V_1 \neq \phi$ and $S \cap V_2 \neq \phi$. Then there exist $x \in S \cap V_1$ and $y \in S \cap V_2$ such that $xy \neq yx$. Now $x + y \in S$ and $x + y \notin C_K(S)$ that is $x +y \in V_1$ or $V_2$, these give contradictions. Hence, $S \cap V_1 = \phi$ or $S \cap V_2 = \phi$. That is $S \subseteq V_2$ or $S \subseteq V_1$. Suppose $S \subseteq V_1$. Then   the vertices of $\Gamma_{S, K}$ belonging to  $S$ are not adjacent to any of the  vertices in $S$. Therefore, if $v \in V_1$ then $vs = sv$ for all $s \in S \setminus C_K(S)$. Thus $v \in Z(S) \subseteq C_K(S)$, a contradiction. Hence the theorem follows.
\end{pf}

\begin{Thm}\label{n-regular}
Let $n$ be any square free odd positive integer. Then there is no non-commutative ring $R$ with subrings $S, K$ and $S \subseteq K$ such that $\Gamma_{S, K}$ is an $n$-regular graph.
\end{Thm}
\begin{pf}
Let $\Gamma_{S, K}$ be an $n$-regular graph, where $S \subseteq K$ are subrings of a non-commutative ring $R$. Suppose $n = p_1p_2\dots p_k$, where $p_i$'s are distinct odd primes. If $r \in S$ is a vertex of $\Gamma_{S, K}$  then 
\[
n = deg(r) = |K| - |C_K(r)| = |C_K(r)|([K : C_K(r)] - 1).
\]
If $|C_K(r)| = 1$ then $n = deg(r) = |K| - 1$. But $|V(\Gamma_{S, K})| = |K| - |C_K(S)| \leq |K| - 1$ which gives $deg(r) \leq |K| - 2$, a contradiction. Hence $|C_K(r)| \neq 1$. Thus $|C_K(r)| = \underset{p_i \in Q}{\prod}p_i$ and $[K : C_K(r)] - 1 = \underset{p_j \in P\setminus Q}{\prod}p_j$, where $Q \subseteq \{p_1, p_2, \dots, p_k\} = P$. So, $|K| = \underset{p_i \in Q}{\prod}p_i(\underset{p_j \in P \setminus Q}{\prod}p_j + 1)$. If $r \in R \setminus S$ then, using similar argument, we have $|S| = \underset{p_i \in T}{\prod}p_i(\underset{p_j \in P\setminus T}{\prod}p_j + 1)$ where $T \subseteq  P$. Since $|S|$ divides $|K|$,   $\underset{p_i \in T - (T \cap Q)}{\prod}p_i(\underset{p_j \in P \setminus T}{\prod}p_j + 1)$ divides $\underset{p_j \in P \setminus Q}{\prod}p_j + 1$, which is not possible. This completes the proof.
\end{pf}

Now putting $S = R$ in Theorems \ref{complete relative non-commuting}-\ref{n-regular}, we have the following corollary.

\begin{Cor}\label{complete non-commuting}
There is no non-commutative ring $R$ with subrings $S, K$ and $S \subseteq K$ such that 
\begin{enumerate}
\item $\Gamma_R$ is a star graph.
\item $\Gamma_R$ is complete bipartite graph.
\item $\Gamma_R$ is an $n$-regular graph where $n$ is any square free odd positive integer.
\end{enumerate}
\end{Cor}

A complete graph is is a graph in which every pair of distinct vertices is adjacent. In the following theorem we show that a complete graph can not be realized by $\Gamma_{S, K}$ for some subrings $S, K$ of $R$. 
\begin{Thm}
There is no non-commutative ring $R$ with subrings $S, K$ where $S \subseteq K$ and $K$ has unity such that $\Gamma_{S, K}$ is complete. In particular, there is no non-commutative ring $R$ with unity such that $\Gamma_{R}$ is complete.
\end{Thm}
\begin{pf}
If $S$ is commutative then $\Gamma_{S, K}$ is an empty graph. Suppose that  $S$ is non-commutative and  $\Gamma_{S, K}$ is complete. Then for $s \in V(\Gamma_{S, K})\cap S$ we have $deg(s) = |V(\Gamma_{S, K})| - 1 = |K| - |C_K(S)| - 1$. By part (a) of Corollary \ref{deg_cor}, we have $|K| - |C_K(s)| = |K| - |C_K(S)| - 1$. This gives $|C_K(S)| = 1$ and $|C_K(s)| = 2$, a contradiction. This completes the proof of the first part.

Particular case follows by putting $S = R$.
\end{pf}

\section{Diameter, girth and dominating set}

In this section, we obtain diameter, girth and dominating set of $\Gamma_{S, K}$. We write $diam(G)$ and $girth(G)$ to denote the diameter and girth of a graph $G$ respectively. Recall that 
$diam(G) = max\{d(x, y) : x, y \in V(G)\}$, where $d(x, y)$ denotes the distance between $x$ and $y$. Also $girth(G)$ is the length of the shortest cycle obtained in $G$.  

\begin{Thm}\label{diam and girth}
Let $R$ be a non-commutative ring and $S, K$ two subrings of $R$.
\bnum
\item If $Z(S) = Z(K) = \{0\}$ then $diam(\Gamma_{S, K}) \leq 3$ and \\$girth(\Gamma_{S, K}) \leq 4$.
\item If $S \subseteq K$ and $Z(S) = \{0\}$ then $diam(\Gamma_{S, K}) = 2$ and \\$girth(\Gamma_{S, K}) = 3$.  
\enum
\end{Thm}
\begin{pf}
(a) Suppose $r \in V(\Gamma_{S, K})$. If $S \subseteq C_K(r)$ or $K \subseteq C_S(r)$ then $r \in C_K(S)$ or $r \in C_S(K)$, a contradiction.
Therefore $S \nsubseteq C_K(r)$ and $K \nsubseteq C_S(r)$. Suppose $r$ and $t$ are vertices of $\Gamma_{S, K}$ such that they are not adjacent.
If $r, t \in K$ then there exist vertices $s_1, s_2 \in S$ such that $rs_1 \neq s_1r$ and $ts_2 \neq s_2t$. If $r$ is adjacent to $s_2$ or $t$ is adjacent to $s_1$, then $d(r, t) = 2$. If $r \in S$ and $t \in K$ then there exist vertices $k \in K$ and $s \in S$ such that $rk \neq kr$ and $ts \neq st$, as $S \nsubseteq C_K(r)$ and $K \nsubseteq C_S(r)$. Suppose $r$ is adjacent to $s$ or $t$ is adjacent to $k$, then $d(r, t) = 2$. If they are not adjacent and $k$ is adjacent to $s$ then $d(r, t) = 3$. If $r$ is not adjacent to $s$, $t$ is not adjacent to $k$ and $k$ is not adjacent to $s$ then $(k + t)$ is adjacent to $r$ and $s$. So $d(r, t) = 3$. If $r, t \in S$ 
then there exist vertices $k_1, k_2 \in K$ such that $rk_1 \neq k_1r$ and $tk_2 \neq k_2t$. If $r$ is adjacent to $k_2$ 
or $t$ is adjacent to $k_1$ then $d(r, t) = 2$. If they are not adjacent then $k_1 + k_2$ is adjacent to $r$ and $t$ and so $d(r, t) = 2$. Hence, $diam(\Gamma_{S, K}) \leq 3$.

Next we suppose that $k \in K, s \in S$ such that $k$ and $s$ are vertices and they are adjacent. So there exist two vertices $s' \in S, k' \in K$ such that $sk' \neq k's$ and $s'k \neq ks'$. If $k$ is adjacent to $k'$ or $s$ is adjacent to $s'$ then $\{s, k, k'\}$ or $\{s, k, s'\}$ is a cycle of length $3$ in $\Gamma_{S, K}$.
Suppose both are not adjacent; and $s', k'$ are adjacent. Then $\{s, k, s', k'\}$ is 
a cycle of length $4$ in $\Gamma_{S, K}$. Suppose $s', k'$ are also not adjacent, then there exists $s + s' \in S$ such that $(s + s')$ is adjacent to $k$ and $k'$. Then $\{s, k, s + s', k'\}$ is   a cycle of length $3$ in $\Gamma_{S, K}$. Hence, $girth(\Gamma_{S, K}) \leq 4$. 

(b) Let $r_1$ and $r_2$ be two vertices of $\Gamma_{S, K}$  such that $r_1r_2 = r_2r_1$. As $r_1$ and $r_2$ are vertices, therefore there exist vertices $s_1, s_2 \in S$ such that $r_1s_1 \neq s_1r_1$ and $r_2s_2 \neq s_2r_2$. If $r_2$ is adjacent to $s_1$ or $r_1$ is adjacent to $s_2$, then $d(r_1, r_2) = 2$. We assume that both are not adjacent, that is $r_1s_2 = s_2r_1$ and $r_2s_1 = s_1r_2$. Then $s_1 + s_2$ is adjacent to $r_1$ and $r_2$, which also gives $d(r_1, r_2) = 2$. Hence, $diam(\Gamma_{S, K}) = 2$.

Next, we suppose that  $r, s \in V(\Gamma_{S, K})$ where $s \in S$ and $r, s$ are adjacent. So, there exist $t_1, t_2 \in V(\Gamma_{S, K})$ such that $rt_1 \neq t_1r$ and $st_2 \neq t_2s$. That is, $r$ and $s$ are adjacent to $t_1$ and $t_2$ respectively. If $r, t_2$ or $s, t_1$ are adjacent then $\{r, s, t_2\}$ or  $\{r, s, t_1\}$ is a cycle of length $3$ in $\Gamma_{S, K}$. Suppose both are not adjacent then it can be seen that $t_1 + t_2$ is adjacent to $r$ and $s$. Therefore, $\{r, s, t_1 + t_2\}$ is a cycle of length $3$ in $\Gamma_{S, K}$.
  Hence, $girth(\Gamma_{S, K})$ is $3$. 
\end{pf}

As a consequence of Theorem \ref{diam and girth}, we have the following corollary.

\begin{Cor}
Let $R$ be a non-commutative ring and $S, K$ two subrings of $R$ such that $S \subseteq K$ and $Z(S) = \{0\}$. Then $\Gamma_{S, K}$ is connected. 
\end{Cor} 
 
Let $G$ be a graph and  $D$  a subset of $V(G)$ such that every vertex not in $D$ is adjacent to at least one member of $D$ then $D$ is called the dominating set for $G$. It is easy to see that for non-commutative subrings $S, K$ of $R$ such that $S \subseteq K$, $V(\Gamma_S)$ is a dominating set for $\Gamma_{S, K}$ and $V(\Gamma_{S, K})$ is a dominating set for $\Gamma_K$ if $|C_K(S)| = 1$. In the next few results we discuss about  dominating sets for $\Gamma_{S, K}$.

\begin{Prop}
Let $S, K$ be two subrings of a non-commutative ring $R$ and
$X \subseteq V(\Gamma_{S, K})$. Then $X$ is a dominating set for $\Gamma_{S, K}$ if and only if $C_S(X) \cup C_K(X) \subseteq X \cup C_S(K) \cup C_K(S)$. 
\end{Prop}
\begin{pf}
Let $X$ be a dominating set for $\Gamma_{S, K}$. Let $r \in V(\Gamma_{S, K})$ such that $r \in C_K(X) \cup C_S(X)$. Also $r \notin C_S(K) \cup C_K(S)$. If $r \notin X$ then there exists an element $x \in X$ such that $rx \neq xr$, a contradiction. If $r \notin V(\Gamma_{S, K})$ such that $r \in C_K(X) \cup C_S(X)$ then $r \in C_S(K) \cup C_K(S)$.

Conversely, we suppose that $C_S(X) \cup C_K(X) \subseteq X \cup C_S(K) \cup C_K(S)$. Let $l \in V(\Gamma_{S, K})$ such that $l \notin X$. Suppose $lx = xl$ for all $x \in X$, that is $X$ is not a dominating set. So, $l \in C_S(X)$ or $l \in C_K(X)$. This gives $l \in X \cup C_S(K) \cup C_K(S)$. Therefore, $l \in X$, a contradiction. Hence, $X$ is a dominating set for $\Gamma_{S, K}$.
\end{pf}

\begin{Prop}
Let $R$ be a non-commutative ring with unity and $S, K$ be two subrings of $R$. Let $A = \{s_1, s_2, \dots, s_m\}$ and $B = \{k_1, k_2, \dots,\\ k_n\}$ be generating sets for $S$ and $K$ respectively. If $A \cap (C_S(K) \cup C_K(S)) = \{s_{c + 1}, \dots, s_m\}$ and $B \cap (C_S(K) \cup C_K(S)) = \{k_{d + 1}, \dots, s_n\}$ then $X = \{s_1, s_2,\dots, s_c, k_1, k_2,\dots, k_d\}$ is a dominating set for $\Gamma_{S, K}$. 
\end{Prop}

\begin{pf}
Clearly $X \subseteq V(\Gamma_{S, K})$. Let $r \in V(\Gamma_{S, K})$ such that $r \notin X$. If $r \in S$ then there exists an element $k \in K$ such that $k = g_ik_1^{\alpha_{1i}}k_2^{\alpha_{2i}}\dots k_p^{\alpha_{pi}}$ where $g_i \in \Z$, $\alpha_{ji} \in {\mathbb{N}} \cup \{0\}$ and $k_j \in B$ such that $rk \neq kr$. Thus $rk_i \neq k_ir$ for some $i, 1 \leq i \leq d$. If $r \in K$ then there exists an element $s \in S$ such that $s = h_js_1^{\alpha_{1j}}s_2^{\alpha_{2j}}\dots k_q^{\alpha_{qj}}$ where $h_j \in \Z$, $\alpha_{ij} \in {\mathbb{N}} \cup \{0\}$ and $s_i \in A$ such that $rs \neq sr$. Thus $rs_j \neq s_jr$ for some $j, 1 \leq j \leq c$. This completes the proof.
\end{pf}

\begin{Prop}
Let $R$ be a non-commutative ring with unity and $S, K$ two subrings of $R$ such that $S \subseteq K$. Let $A = \{s_1, s_2, \dots, s_n\}$ be a generating set for $S$. If $A \cap C_K(S) = \{s_{m + 1}, \dots, s_n\}$ then $B = \{s_1, s_2, \dots, s_m\} \cup \{s_1 + s_{m + 1}, s_1 + s_{m + 2}, \dots, s_1 + s_n\}$ is a dominating set for $\Gamma_{S, K}$. 
\end{Prop} 

\begin{pf}
Clearly $B \subseteq V(\Gamma_{S, K})$. Let $r$ be an element of $V(\Gamma_{S, K})$ such that $r \notin B$. If $r \in S$ then there exists an element $s =  z_i{s_1}^{\alpha_{1i}}{s_2}^{\alpha_{2i}}\dots  {s_d}^{\alpha_{di}}$, where $z_i \in \Z$, $\alpha_{ji} \in \N$ $\cup \{0\}$ and $S_j \in A$ such that $rs \neq sr$. Hence, $rs_i \neq s_ir$ for some $1 \leq i \leq m$ and so $r$ is adjacent to $s_i$. 

If $r \in K \setminus S$ then there exists an element $t = z_i{s_1}^{\alpha_{1i}}{s_2}^{\alpha_{2i}}\dots {s_p}^{\alpha_{pi}}$, where $z_i \in \Z,$ $\alpha_{li} \in \N$ $\cup \{0\}$ and $S_l \in A$ such that $rt \neq tr$. If $rs_i \neq s_ir$ for some $1 \leq i \leq m$ then $r$ is adjacent to $s_i$. Otherwise, $rs_i = s_ir$ for all $1 \leq i \leq m$. Since $r \notin C_K(S)$, there exists $s_l$ for some $ m + 1 \leq l \leq n$ such that $rs_l \neq s_lr$. Hence, $r$ is adjacent to $s_1 + s_l$. This completes the proof.
\end{pf}
We conclude this section by the following result.
\begin{Prop}
Let $S, K$ be two non-commutative subrings of a ring $R$ such that $S \subseteq K$. Then $X = (S + C_K(S)) \setminus C_K(S)$ is a dominating set for $\Gamma_{S, K}$.
\end{Prop}

\begin{pf}
Suppose $r$ is a vertex of $\Gamma_{S, K}$ such that $r \notin X$. So there exists an element $s \in S$ such that $rs \neq sr$. If $s \notin C_K(S)$ then $s \in X$ and $s$ is adjacent to $r$. If $s \in C_K(S)$ then there exists an element $t \in S \setminus C_K(S)$ such that $st \neq ts$. If $rt \neq tr$ then $r$ is adjacent to $t$ and $t \in X$. If $rt = tr$ then $t + s \in (S + C_K(S)) \setminus C_K(S)$ and $r(t + s) \neq (t + s)r$. So $t + s \in X$ and $r$ is adjacent to $t + s$.
This completes the proof.    
\end{pf}

\section{Relation between $\Gamma_{S, K}$ and $\Pr(S, K)$}
If $R_1$ and $R_2$ are two non-commutative rings with centers of equal order such that $\Gamma_{R_1}$ and $\Gamma_{R_2}$ are isomorphic graphs then it is easy to see that their commuting probabilities are same. In this section, we give some more connections between $\Gamma_{S, K}$ and $\Pr(S, K)$, where $S, K$ are subrings of $R$. We begin with the following result. 

\begin{Thm}\label{relation}
Let $S$ and $K$ be two subrings of a non-commutative ring $R$ such that $S \subseteq K$. Then the number of edges of $\Gamma_{S, K}$ is
\begin{center}
$
|E(\Gamma_{S, K})| = |S||K|(1 - \Pr(S, K)) - \frac{|S|^2}{2}(1 - \Pr(S)). 
$
\end{center}
\end{Thm}

\begin{pf}
Let $A = \{(r_1, r_2) \in S \times K : r_1r_2 \neq r_2r_1\}$ and $B = \{(r_1, r_2) \in K \times S : r_1r_2 \neq r_2r_1\}$.   We have 
\begin{align*}
|A| &= |S||K| - |\{(r_1, r_2) \in S \times K : r_1r_2 = r_2r_1\}|\\ 
&= |S||K| - |S||K|\Pr(S, K) = |B|
\end{align*}
and 
\[
|A \cap B| = |\{(a, b) \in S \times S : ab \neq ba\}| = |S|^2 - |S|^2 \Pr(S).
\]
 Hence, the result follows from the fact that  $|E(\Gamma_{S, K})| = |A \cup B|$.
\end{pf}
Putting $S = R$ in Theorem \ref{relation}, we have the following  corollary.
\begin{Cor}\label{relation E&R}
Let $R$ be a non-commutative ring. Then the number of edges of $\Gamma_R$ is

\begin{center}
$
|E(\Gamma_R)| = \frac{|R|^2}{2}(1 - \Pr(R)). 
$
\end{center}
\end{Cor}
Using similar techniques, as in the proof of  \cite[Proposition $3.1$]{GET14}, we also have the following result.
\begin{Thm}\label{relationEandPr}
Let $S$ and $K$ be two subrings of a non-commutative ring $R$ such that $S \nsubseteq K$. Then

$|E(\Gamma_{S, K})| = |S||K|(1 - \Pr(S, K)) + \frac{|S|^2(1 - \Pr(S))}{2} - \frac{|S \cap K|^2(1 - \Pr(S \cap K))}{2}$. 
\end{Thm}
In view of  the above results we have that  lower or upper bounds for $\Pr(R), \Pr(S)$ and $\Pr(S, K)$ will give lower or upper bounds for $|E(\Gamma_R)|, |E(\Gamma_{S, K})|$ and vice-versa. As an example, we have the following lower bound for $|E(\Gamma_{S, K})|$.

\begin{Cor}
Let $S$ and $K$ be two non-commutative subrings of a ring $R$ such that $S \nsubseteq K$. Then 
\[
|E(\Gamma_{S, K})| \geq \frac{3|S|(|K| + |S|/2)}{8} - \frac{|S \cap K|^2}{2}.
\]
\end{Cor}
\begin{pf}
Using  \cite[Theorem $2.4$]{duttabasnetnath} and \cite[Theorem $1$]{machale}, we have $\Pr(S, K) \leq \Pr(S) \leq \frac{5}{8}$. Therefore, by Theorem \ref{relationEandPr}, we have
\begin{align*}
|E(\Gamma_{S, K})| &\geq \frac{3|S||K|}{8} + \frac{3|S|^2}{16} - \frac{|S \cap K|^2}{2} + \frac{|S \cap K|^2\Pr(S \cap K)}{2}\\
&\geq\frac{3|S|( |K| + |S|/2 )}{8} - \frac{|S \cap K|^2}{2}.
\end{align*}
\end{pf}

\begin{Prop}
Let $R$ be a non-commutative ring. Then 
\begin{center}
  $\Pr(R) \geq \frac{2|Z(R)|}{|R|} + \frac{1}{|R|} -  \frac{|Z(R)|^2}{|R|^2} - \frac{|Z(R)|}{|R|^2}$. 
  \end{center}  
\end{Prop}

\begin{pf}
We know that for every graph, the number of edges is at most $\frac{n(n - 1)}{2}$, where $n$ is the number of vertices of the graph. Therefore, $|E(\Gamma_R)| \leq \frac{1}{2}((|R| - |Z(R)|)(|R| - |Z(R)| - 1))$. Hence, using Corollary \ref{relation E&R}, we have the required result. 
\end{pf} 

\begin{Prop}\label{lower bound for E}
Let $S$ and $K$ be two subrings of a non-commutative ring $R$ such that $S \subseteq K$. Then
\[
|E(\Gamma_{S, K})| \geq \frac{1}{2}|S||K| - \frac{1}{4}|S|^2 - \frac{1}{4}|Z(S)||K| - \frac{1}{4}|S||C_K(S)| + \frac{1}{4}|Z(S)||S|.
\]
\end{Prop}

\begin{pf}
Let $S_1 = V(\Gamma_{S, K}) \cap S$ and $S_2 = V(\Gamma_{S, K}) \setminus S$. Then $|S_1| = |S| - |Z(S)|$ and $|S_2| = |K| - |S| - |C_K(S)| + |Z(S)|$. We have
\begin{align*}
|E(\Gamma_{S, K})| &= \underset{r \in V(\Gamma_{S, R})}{\sum}deg(r) = \underset{r \in S_1}{\sum}deg(r) + \underset{r \in S_2}{\sum}deg(r)\\
&= \underset{r \in S_1}{\sum}(|K| - |C_K(r)|) + \underset{r \in S_2}{\sum}(|S| - |C_S(r)|)\\ &\geq |S_1||K| - \frac{|S_1||K|}{2} - |S_2||S| - \frac{|S||S_2|}{2}.
\end{align*}
Now putting the values of $|S_1|$ and $|S_2|$  we have the required result.    
\end{pf}

Putting $S  = R$ in Proposition \ref{lower bound for E} and then using Corollary \ref{relation E&R} we have the   upper bound for $\Pr(R)$.
\begin{Cor}\label{up2}
Let $R$ be a non-commutative ring. Then 
\begin{center}
$\Pr(R) \leq \frac{1}{2} + \frac{1}{2}\frac{|Z(R)|}{|R|}$.
\end{center}
\end{Cor}
Note that the upper bound obtained in Corollary \ref{up2} is slightly better than the upper bound obtained in \cite[Theorem 1]{machale}.

\begin{Prop}
Let $S$ and $K$ be two non-commutative subrings of a ring $R$ such that $S \subseteq K$ and $p$ the smallest prime dividing $|R|$. Then 
\[
|E(\Gamma_{S, K})| \leq |S|(|K| - \frac{3|S|}{16} - p) - |Z(K) \cap S|(|K| - p)
\]
\end{Prop}
\begin{pf}
By \cite[Theorem 2.5]{duttabasnetnath}, we have 
\[
\frac{|Z(K) \cap S|}{|S|} + \frac{p(|S| - |Z(K) \cap S|)}{|S||K|} \leq \Pr(S, K).
\]
Now using this and the fact that $\Pr(S) \leq \frac{5}{8}$  in Theorem \ref{relation}, we get the required result. 
\end{pf}

\begin{Prop}
Let $S$ and $K$ be two non-commutative subrings of a ring $R$ such that $S \subseteq K$. Then 
\[
|E(\Gamma_{S, K})| \geq -\frac{3|S|^2}{16} + \frac{3|S||K|}{8}.
\]
\end{Prop}
\begin{pf}
Using  \cite[Theorem 2.2 ]{duttabasnetnath}, we have that $\Pr(S, K) \leq \Pr(S) \leq \frac{5}{8}$, as $S$ is non-commutative. So, $1 - \Pr(S, K) \geq 1 - \Pr(S) \geq \frac{3}{8}$. Therefore, the result follows from Theorem \ref{relation}. 
\end{pf}

\begin{Prop}
There is no non-commutative ring $R$ with non-commutative subrings $S, K$ such that $S \subseteq K$ and $|C_K(S)| = 1$ satisfying
\[
2|K|\Pr(S, K) - |S|\Pr(S) = -2|K|/|S| + 4/|S| + 2|K| - |S|.
\]
In particular, there is no non-commutative ring $R$ with trivial center having commuting probability $\Pr(R) = 1 - 2/|R| + 4/|R|^2$.
\end{Prop}
\begin{pf}
Suppose there exists a non-commutative ring $R$ with non commutative subrings $S, K$ such that $S \subseteq K$ and
\[
2|K|\Pr(S, K) - |S|\Pr(S) = -2|K|/|S| + 4/|S| + 2|K| - |S|.
\]
If $|C_K(S)| = 1$ then the above equation, in view of Theorem \ref{relation}, gives
 \[
|E(\Gamma_{S, K})| = |K| - |C_K(S)| - 1 = |V(\Gamma_{S, K})| -1. 
\]
With this relation we can easily create a star graph $\Gamma_{S, K}$, which is a contradiction      (by Theorem \ref{star}). This proves the first part of the proposition.

 Second part is obtained by putting  $S =  R$. 
\end{pf}
We conclude this section by the following result.
\begin{Prop}
There is no non-commutative ring $R$ with commuting probability
\[
\Pr(R) = \frac{1}{2} + \frac{|Z(R)|}{|R|} - \frac{|Z(R)|^2}{2|R|^2}.
\]
\end{Prop}

\begin{pf}
Suppose there exists a non-commutative ring $R$ such that $\Pr(R) = \frac{1}{2} + \frac{|Z(R)|}{|R|} - \frac{|Z(R)|^2}{2|R|^2}$. Then $\frac{|R|^2}{2}(1 - \Pr(R)) = \frac{1}{4}(|R| - |Z(R)|)^2 = |E(\Gamma_R)|$. Therefore   $V(\Gamma_R) = R\setminus Z(R)$ can be partitioned equally into two disjoint sets such that each vertex of one set is adjacent to all vertices of the other set. Thus $\Gamma_R$ is a complete bipartite graph, which is not possible (by  part (b) of Corollary \ref{complete non-commuting}). Hence the result follows. 
\end{pf}

\section{Relation between $\Z$-isoclinism and $\Gamma_{S, K}$}

Hall \cite{pH40} introduced the notion of isoclinism between two groups and Lescot \cite{pL95} showed that the commuting probability of two isoclinic finite groups are same. Later on Buckley, MacHale and N$\acute{\rm i}$ sh$\acute{\rm e}$ \cite{BMS} introduced the concept of $\Z$-isoclinism between two rings and showed that the commuting probability of two isoclinic finite rings are same. In \cite{duttabasnetnath}, we introduce the concept of $\Z$-isoclinism between two pairs of rings and show that relative commuting probability remains invariant under $\Z$-isoclinism of pairs of rings.
In this section, we further generalize $\Z$-isoclinism between pairs of rings and find some connections between these pairs and their generalized non-commuting graphs.  

\begin{Def} 
Let $R_1$ and $R_2$ be two rings with subrings $S_1, K_1$ and $S_2, K_2$ respectively such that $S_1 \subseteq K_1$ and $S_2 \subseteq K_2$. A pair of rings $(S_1, K_1)$ is said to be $\Z$-isoclinic to a pair of rings $(S_2, K_2)$ if there exist additive group isomorphisms $\phi : \frac{K_1}{Z(K_1) \cap S_1} \rightarrow \frac{K_2}{Z(K_2) \cap S_2}$ such that $\phi(\frac{S_1}{Z(K_1) \cap S_1}) = \frac{S_2}{Z(K_2) \cap S_2}$ and $\psi : [S_1, K_1] \rightarrow [S_2, K_2]$ such that $\psi([u, v]) = [u', v']$ whenever $u \in S_1, u' \in S_2, v \in K_1, v' \in K_2, \phi(u + (Z(K_1) \cap S_1)) = u' + (Z(K_2) \cap S_2)$ and $\phi(v + (Z(K_1) \cap S_1)) = v' + (Z(K_2) \cap S_2)$. Such pair of mappings $(\phi, \psi)$ is called a generalized $\Z$-isoclinism from $(S_1, K_1)$ to $(S_2, K_2)$.
\end{Def}
We have the following main result of this section.
\begin{Thm} 
Let $R_1$ and $R_2$ be two rings with subrings $S_1, K_1$ and $S_2, K_2$ respectively such that $S_1 \subseteq K_1$ and $S_2 \subseteq K_2$. Let the pairs $(S_1, K_1)$ and $(S_2, K_2)$ are generalized $\Z$-isoclinic. Then $\Gamma_{S_1, K_1} \cong \Gamma_{S_2, K_2}$ if $|Z(K_1) \cap S_1| = |Z(K_2) \cap S_2|$ and $|Z(K_1)| = |Z(K_2)|$.  
\end{Thm}

\begin{pf}
Let $(\phi, \psi)$ be a generalized $\Z$-isoclinism between $(S_1, K_1)$ and $(S_2, K_2)$. Then $|\frac{K_1}{Z(K_1) \cap S_1}| = |\frac{K_2}{Z(K_2) \cap S_2}|,$ $|\frac{S_1}{Z(K_1) \cap S_1}| = |\frac{S_2}{Z(K_2) \cap S_2}|$ and $|[S_1, K_1]| = |[S_2, K_2]|$. Therefore $|S_1| = |S_2|, |\frac{K_1}{Z(K_1)}| = |\frac{K_2}{Z(K_2)}|, |Z(K_1) \setminus S_1| = |Z(K_2) \setminus S_2|$ and $|S_1 \setminus Z(K_1)| = |S_2 \setminus Z(K_2)|$, as $|Z(K_1) \cap S_1| = |Z(K_2) \cap S_2|$ and $|Z(K_1)| = |Z(K_2)|$. Also, by second isomorphism theorem, we have $\frac{S_1}{S_1 \cap Z(K_1)} \cong \frac{S_1 + Z(K_1)}{Z(K_1)}$ (additive group isomorphism). Let $\{s_1, s_2, \dots, s_k\}$ be a transversal for $\frac{S_1 + Z(K_1)}{Z(K_1)}$. Then the set $\{s_1, s_2, \dots, s_k\}$ can be extended to a transversal for $\frac{K_1}{Z(K_1)}$. Let $\{s_1, s_2, \dots, s_k, r_{k + 1}, \dots, r_n\}$ be a transversal for $\frac{K_1}{Z(K_1)}$. Similarly, we can find a transversal $\{s'_1, s'_2, \dots, s'_k, r'_{k + 1}, \dots, r'_n\}$ for $\frac{K_2}{Z(K_2)}$ such that $\{s'_1, s'_2, \dots, s'_k\}$ is a transversal for $\frac{S_2 + Z(K_2)}{Z(K_2)} \cong \frac{S_2}{S_2 \cap Z(K_2)}$.

Let $\phi$ be defined as $\phi(s_i + Z(K_1)) = s'_i + Z(K_2)$, $\phi(r_j + Z(K_1)) = r'_j + Z(K_2)$ for $1 \leq i \leq k, k + 1 \leq j \leq n$ and let the one-to-one correspondence $\theta : Z(K_1) \rightarrow Z(K_2)$ maps elements of $S_1$ to $S_2$. Then $|C_{K_1}(S_1)| = |C_{K_2}(S_2)|$. Let us define a map $\alpha : K_1 \rightarrow K_2$ such that $\alpha(s_i + z) = s'_i + \theta(z)$, $\alpha(r_j + z) = r'_j + \theta(z)$ for $1 \leq i \leq k, k + 1 \leq j \leq n$ and $z \in Z(K_1)$. Then $\alpha$ is a bijection. This shows that $\alpha$ is also a bijection from $K_1 \setminus C_{K_1}(S_1)$ to $K_2 \setminus C_{K_2}(S_2)$. Suppose $r_1, r_2$ are adjacent in $\Gamma_{S_1, K_1}$. Then $r_1 \in S_1$ or $r_2 \in S_1$, say $r_1 \in S_1$. So, $[r_1, r_2] \neq 0$, this gives $[s_i + z, r +z_1] \neq 0$, where $r_1 = s_i + z, r_2 = r +z_1$ for some $z, z_1 \in Z(K_1)$, $r \in \{s_1, s_2, \dots, s_k, r_{k + 1}, \dots, r_n\}$ and $1 \leq i \leq k$. Thus $[s'_i + \theta(z), r + \theta(z_1)] \neq 0$, where $\theta(z), \theta(z_1) \in Z(K_2)$. Hence $[\alpha(s_i + z), \alpha(r + z_1)] \neq 0$, that is $\alpha(r_1)$ and $\alpha(r_2)$ are adjacent. This completes the proof of the theorem. 
\end{pf}

\begin{Thm}
Let $S_1$ and $S_2$ be two subrings of a non-commutative ring $R$ such that $\Gamma_{S_1, R} \cong \Gamma_{S_2, R}$. Then $\Gamma_{S_1} \cong \Gamma_{S_2}$.
\end{Thm}
\begin{pf}
Let $\phi$ be an isomorphism between $\Gamma_{S_1, R}$ and $\Gamma_{S_2, R}$. Suppose that there exists an element $s \in V(\Gamma_{S_1, R}) \cap S_1$ such that $\phi(s) \in V(\Gamma_{S_2, R}) \cap (R \setminus S_2)$. We have that $deg(s) = deg(\phi(s))$. This gives $|R| - |C_R(S_1)| = |S_2| - |C_{S_2}(\phi(s))| < |S_2|$ and so $\frac{|R|}{2} < |S_2|$. Thus $|S_2| = |R|$, a contradiction. Therefore $\phi$ is a bijection between   $V(\Gamma_{S_1, R}) \cap S_1$ and $V(\Gamma_{S_2, R}) \cap S_2$. This completes the proof.
\end{pf}

We conclude the paper with the following corollary
\begin{Cor}
Let $S_1$ and $S_2$ be two subrings of a non-commutative ring $R$ such that $(S_1, R)$ is generalized $\Z$-isoclinic to $(S_2, R)$. Then $\Gamma_{S_1} \cong \Gamma_{S_2}$ if $|Z(R) \cap S_1| = |Z(R) \cap S_2|$.
\end{Cor}


\end{document}